\documentclass[11pt]{amsart}
\usepackage{amssymb,pb-diagram}
\usepackage[?]{amsrefs}

\newtheorem{Thm}{Theorem}[section]

\newtheorem{Cor}[Thm]{Corollary}
\newtheorem{Lem}[Thm]{Lemma}
\newtheorem{Prop}[Thm]{Proposition}

\theoremstyle{definition}

\theoremstyle{remark}

\def\ldots{\mathinner{\ldotp\ldotp\ldotp}}
\def\cdots{\mathinner{\cdotp\cdotp\cdotp}}

\def \cal{\mathcal}

\def\tfrac{\textstyle\frac}
\def \diam{\text{diam }}
\def\Lip{\text{Lip}}

\def\tfrac#1#2{{\textstyle\frac#1#2}}

\def\car{\textrm{\bbold 1}}
\font\bbold=bbold12

\newcommand{\R}{\ensuremath{\mathbb{R}}}
\newcommand{\N}{\ensuremath{\mathbb{N}}}

\begin{document}

\title{Best constants for Lipschitz embeddings of metric spaces into $c_0$}
\author{N. J. Kalton}
\address{Department of Mathematics \\
University of Missouri-Columbia \\
Columbia, MO 65211 }

\email{nigel@math.missouri.edu}

\author{G. Lancien}

\address{Universit\'e de Franche-Comt\'e, Laboratoire de Math\'ematiques UMR 6623,
16 route de Gray, 25030 Besan\c con Cedex, FRANCE.}
\email{gilles.lancien@univ-fcomte.fr}


\thanks{The first author acknowledges support from NSF grant DMS-0555670}

\begin{abstract}  We answer a question of Aharoni by showing that every
separable metric space can be Lipschitz 2-embedded into $c_0$ and this result
is sharp; this improves earlier estimates of Aharoni, Assouad and Pelant. We
use our methods to examine the best constant for Lipschitz embeddings of the
classical $\ell_p-$spaces into $c_0$ and give other applications. We prove that
if a Banach space embeds almost isometrically into $c_0$, then it embeds
linearly almost isometrically into $c_0$. We also study Lipschitz embeddings
into $c_0^+$.
\end{abstract}

\maketitle

\section{Introduction}  In 1974, Aharoni \cite{Aharoni1974} proved that every
separable metric space $(M,d)$ is Lipschitz isomorphic to a subset of the
Banach space $c_0.$  Thus, for some constant $K$, there is a map $f:M\to c_0$
which satisfies the inequality
$$ d(x,y)\le \|f(x)-f(y)\|\le Kd(x,y) \qquad x,y\in M.$$ Aharoni proved this
result with $K=6+\epsilon$ where $\epsilon>0$, so that every separable metric
space  $(6+\epsilon)$-embeds into $c_0.$  He also noted that if one takes $M$
to be the Banach space $\ell_1$ one cannot have $K<2.$  In fact the map defined
by Aharoni took values in the positive cone $c_0^+$ of $c_0.$  Later Assouad
\cite{Assouad1978} refined Aharoni's result by showing that every separable
metric space  $(3+\epsilon)$-embeds into $c_0^+$ (see
\cite{BenyaminiLindenstrauss2000} p. 176ff).   A further improvement was
obtained by Pelant in 1994 \cite{Pelant1994} who showed that every separable
metric space  $3$-embeds into $c_0^+$ and that this result is sharp in the
sense that $\ell_1$ cannot be $(3-\epsilon)$-embedded into $c_0^+$ (see also
\cite{Aharoni1978} for the lower bound).

These results leave open the question of the best constant for Lipschitz
embeddings into $c_0.$   Note that $c_0$ can only be $2$-embedded into $c_0^+$.
The main result of this paper is that every separable metric space $2$-embeds
into $c_0$ and this is sharp by Aharoni's remark above. To prove this result,
for $1<\lambda\le 2$ we establish a criterion $\Pi(\lambda)$ which is
sufficient to imply that a separable metric space $\lambda$-embeds into $c_0$
(and the converse statement is almost true). This criterion enables us to
establish sharp results concerning the embedding of $\ell_p$-spaces into $c_0$:
thus $\ell_p$ $2^{1/p}$-embeds into $c_0$ if $1\le p<\infty$ and the constant
is best possible. Using a previous work of the first author and D. Werner
\cite{KaltonWerner1995}, we also show that a Banach space which embeds almost
isometrically into $c_0$ embeds linearly almost isometrically into $c_0$.

The same techniques can be applied to embeddings into $c_0^+$. Here we show
that $\ell_p$ $(2^p+1)^{1/p}$-embeds into $c_0^+$ and $\ell_p^+$
$3^{1/p}$-embeds into $c_0^+$ and in each case the result is best possible.

We conclude the paper by proving that every separable ultrametric space embeds
isometrically into $c_0^+$ and the infinite branching tree embeds isometrically
into $c_0.$

\section{Lipschitz embeddings into $c_0$}

Let $(M,d)$ be a metric space and let $A$ and $B$ be non-empty subsets of $M.$

\noindent We define
$$ \delta(A,B)=\inf_{a\in A,b\in B}d(a,b)$$
and
$$ D(A,B)=\sup_{a\in A,b\in B}d(a,b).$$
In this paper all metric balls are closed with strictly positive radii.  If
$f:(M_1,d_1)\to (M_2,d_2)$ is a Lipschitz map between metric spaces we write
$\Lip(f)$ for the Lipschitz constant of $f$, i.e. the least constant $K$ such
that $d_2(f(x),f(y))\le Kd_1(x,y)$ for $x,y\in M_1.$

\begin{Lem}\label{Lipfunction}  Let $(M,d)$ be a metric space and suppose that $A,B$
and $C$ are non-empty  subsets of $M.$  Then for $\epsilon>0$, there exists a
Lipschitz function $f:M\to\mathbb R$ with $\Lip(f)\le 1$ such that
\newline
(i) $|f(x)|\le \epsilon \qquad x\in C$\newline and\newline (ii) $|f(x)-f(y)|=
\theta=\min\,\big(\delta(A,B),\delta(A,C)+\delta(B,C)+2\epsilon\big) \qquad
x\in A,\ y\in B.$
\end{Lem}

\begin{proof} Let us augment $M$ by adding an extra point $0$; let $M^*=M\cup\{0\}$.
We define:

$$ d^*(x,y)=\begin{cases} \min\,\big(d(x,y),d(x,C)+d(y,C)+2\epsilon\big) \qquad x,y\in M
\\ d(x,C)+\epsilon \qquad x\in M,\ y=0\\ d(y,C)+\epsilon \qquad x=0,\ y\in M
\\ 0 \qquad x=y=0.\end{cases}$$
One can easily check that $d^*$ is a metric on $M^*$. We can pick $s,t$ in $\R$
such that:
$$-(\delta(B,C)+\epsilon)\leq s\leq 0\leq t\leq \delta(A,C)+\epsilon\ \ {\rm
and}\ \ t-s=\theta.$$ Then we define $g: A\cup B\cup\{0\}\to \R$ by $g=t$ on A,
$g=s$ on $B$ and $g(0)=0$. This function is 1-Lipschitz for $d^*$ and can be
extended into a 1-Lipschitz function $f^*$ on $(M^*,d^*)$. Let $f$ be the
restriction of $f^*$ to $M$. Then $f$ satisfies the conditions of the
Lemma.\end{proof}

For $\lambda>1,$ we say that a metric space $(M,d)$ has property $\Pi(\lambda)$
given any $\mu>\lambda$ there exists $\nu>\mu$ such that if $B_1$ and $B_2$ are
two metric balls of radii $r_1,r_2$ respectively then there are finitely many
sets $(U_j)_{j=1}^N,(V_j)_{j=1}^N$ such that:
$$  \lambda\delta(U_j,V_j) \ge \nu(r_1+r_2) \qquad 1\le j\le N$$ and
$$ \{(x,y)\in B_1\times B_2:\ d(x,y)>\mu(r_1+r_2)\} \subset \bigcup_{j=1}^N(U_j\times V_j).$$

\noindent In this definition the sets $U_j,V_j$ are allowed to be repeated. It
is clearly possible, without loss of generality, to assume they are closed.  We
can also (altering the value of $\nu$) assume that they are open.

\begin{Lem}\label{everymetric}  Every metric space has property $\Pi(2).$\end{Lem}

\begin{proof}  For $\mu>2,$ let
$$ U=B_1\cap\{x:\ \exists y\in B_2,\ d(x,y)> \mu(r_1+r_2)\}$$ and
$$ V=B_2\cap\{y:\ \exists x\in B_1,\ d(x,y)>\mu (r_1+r_2)\}.$$
Then $$ \{(x,y)\in B_1\times B_2:\ d(x,y)>\mu(r_1+r_2)\} \subset U\times V.$$
Suppose $x\in U,\ y\in V.$  Let us assume, without loss of generality, that
$r_1\le r_2.$  Then there exists $x'\in U$ with $d(x',y)>\mu(r_1+r_2).$  Hence
$$ d(x,y)>\mu(r_1+r_2)-d(x,x')\ge \mu(r_1+r_2)-2r_1\ge (\mu-1)(r_1+r_2).$$  Therefore we can take $\nu=2\mu-2>\mu.$
\end{proof}

We say that a metric space is {\it locally compact} (respectively, {\it locally
finite}) if all its metric balls are relatively compact (respectively, finite).

\begin{Lem}\label{everylocally} For any $\lambda>1$,
 every locally compact metric space has property
$\Pi(\lambda)$.\end{Lem}

\begin{proof} Let $\mu>\lambda>1$ and $B_1$, $B_2$ be two balls of a locally compact metric space $(M,d)$,
with respective radii $r_1$ and $r_2$. Pick $\nu$ such that
$\mu<\nu<\lambda\mu$. We denote $\Delta=\{(x,y)\in B_1\times B_2:\
d(x,y)>\mu(r_1+r_2)\}.$ Let $\epsilon>0$. Since $M$ is locally compact, there
are finitely many points $(x_j,y_j)_{j=1}^N$ in $\Delta$ such that
$$\Delta \subset \bigcup_{j=1}^N(U_j\times V_j),\ \ {\rm where}\ \
U_j=B(x_j,\epsilon)\ {\rm and}\ V_j=B(y_j,\epsilon).$$ Then, for all $1\leq
j\leq N$, $\lambda\delta(U_j,V_j)>\lambda\mu(r_1+r_2)-2\lambda\epsilon> \nu
(r_1+r_2)$, if $\epsilon$ was chosen small enough, namely
$\epsilon<(2\lambda)^{-1}(\lambda\mu-\nu)(r_1+r_2)$.
\end{proof}

\begin{Prop}\label{lambdaembedding}  Let $\lambda_0\geq 1$. If a metric space
$(M,d)$ $\lambda_0$-embeds into $c_0$ then it has property $\Pi(\lambda)$ for
every $\lambda>\lambda_0$.
\end{Prop}

\begin{proof}  Suppose $\mu>\lambda.$ Let $B_1,B_2$ be metric balls of radii $r_1,r_2$ and centers $a_1,a_2.$
Let $\Delta=\{(x,y)\in B_1\times B_2:\ d(x,y)>\mu(r_1+r_2)\}.$ Let $f:M\to c_0$
be an embedding such that
$$ d(x,y)\le \|f(x)-f(y)\|\le \lambda_0 d(x,y) \qquad x,y\in M.$$
Suppose $f(x)=(f_i(x))_{i=1}^{\infty}.$  Then there exists $n$ so that
$$|f_i(a_1)-f_i(a_2)|<(\mu-\lambda)(r_1+r_2) \qquad i\ge n+1.$$
Thus if $(x,y)\in \Delta$ we have
$$ |f_i(x)-f_i(y)|<(\mu-\lambda)(r_1+r_2)+\lambda_0 r_1+\lambda_0 r_2<d(x,y), \qquad i\ge n+1.$$  Hence
$$ d(x,y)\le \max_{1\le i\le n}|f_i(x)-f_i(y)| \qquad (x,y)\in\Delta.$$

Choose $\epsilon>0$ so that $\lambda(\mu-\epsilon)>\lambda_0\mu.$ By a
compactness argument we can find coverings $(W_k)_{k=1}^m$ of $B_1$ and
$(W'_k)_{k=1}^{m'}$ of $B_2$ such that  we have
$$ |f_i(x)-f_i(x')|\le \tfrac12\epsilon(r_1+r_2) \qquad x,x'\in W_k,\ 1\le i\le n,\ 1\le k\le m,$$ and
$$ |f_i(x)-f_i(x')|\le \tfrac12\epsilon(r_1+r_2) \qquad x,x'\in W'_k,\ 1\le i\le n,\ 1\le k\le m'.$$
Let
$$ \mathcal S=\{(k,k')\ 1\le k\le m,\ 1\le k'\le m':\ W_k\times W'_{k'}\cap\Delta \neq \emptyset\}$$
and then we define $(U_j)_{j=1}^N,\ (V_j)_{j=1}^N$ in such a way that
$(U_j\times V_j)_{j=1}^N$ is an enumeration of $(W_k\times
W_{k'})_{(k,k')\in\mathcal S}.$ Clearly $ \Delta\subset \cup_{j=1}^NU_j\times
V_j.$  Now suppose $x\in U_j,\ y\in V_j.$ Then there exist $x'\in U_j,\ y'\in
V_j$ so that $d(x',y')>\mu(r_1+r_2).$ Thus there exists $i,\ 1\le i\le n$ so
that $|f_i(x')-f_i(y')|>\mu(r_1+r_2).$  However
$$ |f_i(x)-f_i(y)|\ge |f_i(x')-f_i(y')|-\epsilon(r_1+r_2)>(\mu-\epsilon)(r_1+r_2).$$
Hence
$$ \delta(U_j,V_j)\ge \frac{(\mu-\epsilon)}{\lambda_0}(r_1+r_2).$$
Thus we can take $\nu=\lambda\lambda_0^{-1}(\mu-\epsilon)>\mu.$
\end{proof}

We next observe that the definition of $\Pi(\lambda)$ implies a formally
stronger conclusion.

\begin{Lem}\label{soupup}  Let $(M,d)$ be a metric space with property $\Pi(\lambda).$
Then for every $\mu>\lambda$ there is a constant $\nu>\mu$ so that if  $B_1$
and $B_2$ are two metric balls of radii $r_1,r_2$ respectively then there are
finitely many sets $(U_j)_{j=1}^N,(V_j)_{j=1}^N$ such that if $(x,y)\in
B_1\times B_2$ and $d(x,y)>\mu(r_1+r_2)$ then there exists $1\le j\le N$ so
that $x\in U_j,\ y\in V_j$ and:
$$  \lambda \mu\delta(U_j,V_j) \ge \nu d(x,y).$$ \end{Lem}

\begin{proof}  By the definition of $\Pi(\lambda)$ there exists $\nu'>\lambda$
so that such that if $B_1$ and $B_2$ are two metric balls of radii $r_1,r_2$
respectively then there are finitely many sets $(U_j)_{j=1}^N,(V_j)_{j=1}^N$
such that:
$$  \lambda\delta(U_j,V_j) \ge \nu'(r_1+r_2) \qquad 1\le j\le N$$ and
$$ \{(x,y)\in B_1\times B_2:\ d(x,y)>\mu(r_1+r_2)\} \subset
\bigcup_{j=1}^N(U_j\times V_j).$$

Suppose $\mu<\nu<\nu'$ and let $\epsilon>0$ be chosen so that
$(1+\epsilon)\nu=\nu'.$  Let $B_1,B_2$ be a pair of metric balls of radii
$r_1,r_2>0$. Let $D=D(B_1,B_2)$ and let $m$ be the greatest integer such that
$(1+\epsilon)^m\mu(r_1+r_2)\le D.$  We define $B^{(k)}_1$ for $0\le k\le m$ to
be the ball with the same center as $B_1$ and radius $(1+\epsilon)^kr_1$;
similarly $B^{(k)}_2$ for $0\le k\le m$ is the ball with the same center as
$B_2$ and radius $(1+\epsilon)^kr_2.$ For each $0\le k\le m$ we may determine
sets $U_{kl},V_{kl}$ for $1\le l\le N_k$ so that
$$ \lambda \delta(U_{kl},V_{kl}) \ge \nu'(1+\epsilon)^k(r_1+r_2)$$ and
$$ \{(x,y)\in B^{(k)}_1\times B^{(k)}_2:\ d(x,y)>\mu(1+\epsilon)^k(r_1+r_2)\}
\subset \bigcup_{l=1}^{N_k}(U_{kl}\times V_{kl}).$$ Now if $x\in B_1,y\in B_2$
with $d(x,y)>\mu(r_1+r_2)$ we may choose $0\le k\le m$ so that
$$ (1+\epsilon)^k\mu(r_1+r_2)<d(x,y)\le (1+\epsilon)^{k+1}\mu(r_1+r_2).$$
Then for a suitable $1\le l\le N_k$ we have $x\in U_{kl},\ y\in V_{kl}$ and
$$ \lambda\mu\delta(U_{kl},V_{kl})\ge \nu'(1+\epsilon)^k\mu(r_1+r_2)
\ge \frac{\nu'}{1+\epsilon}d(x,y)=\nu d(x,y).$$ Relabeling the sets
$(U_{kl},V_{kl})_{l\le N_k,\ 0\le k\le m}$ gives the conclusion.\end{proof}

\begin{Lem}\label{construction}  Suppose $(M,d)$ has property $\Pi(\lambda)$.
Suppose $0<\alpha<\beta$.  Let $F,G$ be finite subsets of $M$ and let
$\Delta(F,G,\alpha,\beta)$ be the set of $(x,y)\in M\times M$ such that
$$ \lambda (d(x,G)+d(y,G))+\alpha\le d(x,y)
<\lambda(d(x,F)+d(y,F))+\beta.$$ Then there is a finite set $\mathcal
F=\mathcal F(F,G,\alpha,\beta)$ of functions $f:M\to\mathbb R$ with $\Lip(f)\le
\lambda$ such that
$$ |f(x)|\le \lambda\beta \qquad x\in F$$ and
$$ d(x,y)< \max_{f\in \mathcal F}|f(x)-f(y)| \qquad (x,y)\in \Delta(F,G,\alpha,\beta).$$
\end{Lem}

\begin{proof}  Let $R$ be the diameter of $G$.  Then for $(x,y)\in \Delta(F,G,\alpha,\beta)$ we have
$$ \lambda (d(x,y)-R)+\alpha\le d(x,y)$$ so that
$$ (\lambda-1)d(x,y)< \lambda R.$$
Hence
$$ d(x,G)+d(y,G)< \frac{R}{\lambda-1}.$$
We next let
$$ \mu=\lambda +\frac{(\lambda-1)\alpha}{2R}$$ and choose $\nu=\nu(\mu)$ according
to the conclusion of Lemma \ref{soupup}.

\noindent We now fix $\epsilon>0$ so that $4\mu\epsilon<\alpha.$

Let $E=\{x:\ d(x,G)<(\lambda-1)^{-1}R\}.$ Since $E$ is metrically bounded and
$F\cup G$ is finite we can partition $E$ into finitely many subsets
$(E_1,\ldots,E_m)$ so that for each $z\in F\cup G$ we have:
$$ |d(x,z)-d(x',z)|\le \epsilon \qquad x,x'\in E_j,\ 1\le j\le m.$$
Since $G$ is finite, for each $j$ there exist $z_j\in G$ and $r_j\geq 0$ so
that
$$ \inf_{x\in E_j}d(x,z_j)=\inf_{x\in E_j}d(x,G)=r_j.$$
Thus $E_j$ is contained in a ball $B_j$ centered at $z_j$ with radius
$r_j+\epsilon.$

\noindent Now for each pair $(j,k)$ we can find finitely many pairs of sets
$(U_{jkl},V_{jkl})_{l=1}^{N_{jk}}$ such that for every $(x,y)\in E_j\times E_k$
with $ d(x,y)>\mu(r_j+r_k+2\epsilon)$ there exists $1\le l\le N_{jk}$ with
$x\in U_{jkl},y\in V_{jkl}$ and
$$ \lambda\mu\delta(U_{jkl},V_{jkl})\ge \nu d(x,y).$$
We may as well assume that $U_{jkl}\subset E_j$ and $V_{jkl}\subset E_k$.

Then we apply Lemma \ref{Lipfunction} to construct Lipschitz functions
$f_{jkl}:M\to\mathbb R$ where $1\le j,k\le m,\ 1\le l\le N_{jk}$ such that
$\Lip(f_{jkl})\le \lambda,$
$$ |f_{jkl}(x)|\le \lambda \beta \qquad x\in F$$ and
$$ |f_{jkl}(x)-f_{jkl}(y)|\ge \lambda\theta_{jkl} \qquad x\in U_{jkl},\ y\in V_{jkl}$$ where
$$ \theta_{jkl}=\min\,\big(\delta(U_{jkl},V_{jkl}),\delta(U_{jkl},F)+
\delta(V_{jkl},F)+2\beta\big).$$ Now let us suppose $(x,y)\in
\Delta(F,G,\alpha,\beta).$  Then there exists $(j,k)$ so that $x\in E_j,\ y\in
E_k.$  Note that
\begin{align*} d(x,y)&\ge \lambda(d(x,G)+d(y,G))+\alpha\\
&\ge \lambda (r_j+r_k)+\alpha\\
&= \mu (r_j+r_k+2\epsilon)+\alpha-2\mu\epsilon-(\mu-\lambda)(r_j+r_k)\\
&\ge \mu(r_j+r_k+2\epsilon)+\alpha -2\mu\epsilon-(\mu-\lambda)(\lambda-1)^{-1}R\\
 &>\mu(r_j+r_k+2\epsilon).\end{align*}
Thus  there exists $1\le l\le N_{jk}$ so that $x\in U_{jkl},y\in V_{jkl}$ and
$$ \lambda\delta(U_{jkl},V_{jkl})\ge \frac{\nu}{\mu} d(x,y)>d(x,y).$$
On the other hand, $\epsilon<\alpha/2<\beta/2$. So
\begin{align*} \lambda(\delta(U_{jkl},F)+\delta(V_{jkl},F)+2\beta) &\ge \lambda(d(x,F)+d(y,F)+2\beta-2\epsilon)\\
&>\lambda (d(x,F)+d(y,F)+\beta)\\
&> d(x,y)+(\lambda-1)\beta.\end{align*}
  Hence
$$ |f_{jkl}(x)-f_{jkl}(y)|\ge \lambda \theta_{jkl}>d(x,y).$$  Thus we can take for $\mathcal F$
the collection of all functions $f_{jkl}$ for $1\le j,k\le m$ and $1\le l\le
N_{jk}.$\end{proof}

We now state our main result.

\begin{Thm}\label{Pi}  If a separable metric space $(M,d)$ has property
$\Pi(\lambda)$ for $\lambda>1$, then there is a Lipschitz embedding $f:M\to
c_0$ with
$$ d(x,y)<\|f(x)-f(y)\|\le \lambda d(x,y) \qquad x,y\in M,\ x\neq y.$$
\end{Thm}

\begin{proof}  Let $(u_n)_{n=1}^{\infty}$ be a countable dense set of
distinct points of $M.$

\noindent Denote $F_k=\{u_1,\ldots,u_{k}\}$ for $n\ge 1.$ Let
$(\epsilon_n)_{n=1}^{\infty}$ be a strictly decreasing sequence with
$\lim_{n\to\infty}\epsilon_n=0$.

Using Lemma \ref{construction} we can find an increasing sequence of integers
$(n_k)_{k=0}^{\infty}$ (with $n_0=0$) and a sequence $(f_j)_{j=1}^{\infty}$ of
Lipschitz functions $f_j:M\to\mathbb R$ with $\Lip(f_j)\le \lambda$ so that
$$ |f_j(x)|\le \lambda\epsilon_k \qquad x\in F_k, \ n_{k-1}<j\le n_k$$ and
if
\begin{equation}\label{selectk} \lambda (d(x,F_{k+1})+d(y,F_{k+1}))+
\epsilon_{k+1}\le d(x,y)<\lambda (d(x,F_k)+d(y,F_k))+\epsilon_k
\end{equation}
then
$$ \max_{n_{k-1}<j\le n_k}|f_j(x)-f_j(y)|>d(x,y).$$

Define the map $f:M\to\ell_{\infty}$ by $f(x)=(f_j(x))_{j=1}^{\infty}.$  Then
$\Lip(f)\le \lambda$ and since $f$ maps each $u_j$ into $c_0$, $f(M)\subset
c_0.$  Furthermore if $x\neq y$ the sequence
$$ \sigma_k=\lambda(d(x,F_k)+d(y,F_k))+\epsilon_k$$ is decreasing with $\sigma_1>d(x,y)$ and $\lim_{k\to\infty}\sigma_k=0.$  Hence there is exactly one choice of $k$ so that \eqref{selectk} holds and thus
$\|f(x)-f(y)\|>d(x,y).$\end{proof}

As a corollary, we obtain the following improvement of Aharoni's theorem.

\begin{Thm}\label{Assouad} For every separable metric space $(M,d)$ there is
a Lipschitz embedding $f:M\to c_0$ so that
$$ d(x,y)<\|f(x)-f(y)\|\le 2 d(x,y) \qquad x,y\in M,\ x\neq y.$$\end{Thm}

\begin{proof} Combine Lemma \ref{everymetric} and Theorem \ref{Pi}.\end{proof}

\medskip\noindent {\it Remark.} It follows from Proposition 3 in Aharoni's
original paper \cite{Aharoni1974} that the above statement is optimal.

\begin{Thm}\label{embedlocal} For every locally compact metric space $(M,d)$
and every $\lambda>1$, $(M,d)$ $\lambda$-embeds into $c_0$. This result is best
possible.
\end{Thm}

\begin{proof} The existence of the embedding follows immediately
from the combination of Lemma \ref{everylocally} and Theorem \ref{Pi}. The
optimality of the statement follows from Proposition 3.2 in \cite{Pelant1994},
where J. Pelant proved that $[0,1]^\N$ equipped with the distance
$d((x_n),(y_n))=\sum 2^{-n}|x_n-y_n|$ cannot be isometrically embedded into
$c_0$. To complete the the picture we shall now give a locally finite
counterexample.

Let $(e_n)_{n=0}^\infty$ be the canonical basis of $\ell_1$ and consider the
following locally finite metric subspace of $\ell_1$:
$M=\{0,e_0\}\cup\{ne_n,e_0+ne_n;\ n\ge 1\}$. Assume that $f=(f_k)_{k=1}^\infty$
is an isometry from $M$ into $c_0$ such that $f(0)=0$. Then for all $n\neq m$
in $\N$, there exists $k=k_{n,m}\geq 1$ such that
$|f_k(e_0+ne_n)-f_k(me_m)|=1+n+m$. Since $f_k(0)=0$, we obtain that there is
$\varepsilon=\varepsilon_{n,m}\in \{-1,1\}$ such that
$f_k(e_0+ne_n)=\varepsilon(1+n)$ and $f_k(me_m)=-\varepsilon m$. Therefore
$f_k(e_0)=\varepsilon$ and $f_k(ne_n)=\varepsilon n$.

\noindent Since $f(e_0)\in c_0$, there exists an integer $K$ such that for all
positive integers $n\neq m$, $k_{n,m}\leq K$. Hence, if $\alpha(k,n)$ is the
signum of $f_k(ne_n)$, we have that there exists $k\leq K$ so that
$\alpha(k,n)\neq \alpha(k,m)$, whenever $1\leq n<m$. But on the other hand,
there is clearly an infinite subset $A$ of $\N$ such that for every $k\leq K$
and every $n,m \in A$, $\alpha(k,n)=\alpha(k,m)$. This is a contradiction.
\end{proof}

\section{Embeddings of classical Banach spaces}

In this section we will consider the best constants for embedding certain
classical spaces into $c_0.$  We start by establishing a lower bound condition,
using the Borsuk-Ulam theorem.

\begin{Prop}\label{lower} Suppose $X$ is a Banach space and that $f:X\to c_0$
is a Lipschitz embedding with constant $\lambda_0.$  Then for any $u\in X$ with
$\|u\|=1$ and any infinite-dimensional subspace $Y$ of $X$ we have
$$ \inf_{\substack{y\in Y\\ \|y\|=1}}\|u+ y\|\le \lambda_0.$$
\end{Prop}

\begin{proof}  It follows from Lemma \ref{lambdaembedding} that $X$ has
property $\Pi(\lambda)$ for any $\lambda>\lambda_0$. Let us consider
$B_1=-u+B_X$ and $B_2=u+B_X$, where $B_X$ denotes the closed unit ball of $X$.
Suppose $\mu>\lambda_0$ and select $\mu>\lambda>\lambda_0.$ Then, for some
$\nu>\mu$, we can find finitely many sets $(U_j,V_j)_{j=1}^N$  (which we can
assume to be closed)  verifying:
$$\lambda \delta(U_j,V_j)\ge 2\nu$$ and
$$ \{(x,y)\in B_1\times B_2:\ \|x-y\|>2\mu\}\subset \cup_{j=1}^NU_j\times V_j.$$
Now let $E$ be any subspace of $X$ of dimension greater than $N$ and let
$$ A_j=\{e\in E:\ \|e\|=1,\ (-u+e,u-e)\in {U}_j\times {V}_j\}.$$  Thus the sets $A_j$ are
all closed subsets of the unit sphere $S_E$ of $E$. Assume that for any $e\in
S_E$, $\|u-e\|>\mu$. Then $A_1\cup\cdots\cup A_N=S_E.$ We now use a classical
corollary of the Borsuk-Ulam theorem which is in fact due to Lyusternik and
Shnirelman \cite{LyusternikShnirelman1930} and predates Borsuk's work (see
\cite{Matousek2003} p. 23).  This gives the existence of $e$ in $S_E$ and
$k\leq N$ such that $e$ and $-e$ belong to $A_k$, i.e. $-u\pm e\in U_k$ and
$u\pm e\in V_k.$ This implies that $\delta (U_k,V_k)\le2$ and hence $\lambda\ge
\nu>\mu$ which is a contradiction. Thus there exists $e\in S_E$ with $\|u-
e\|\le \mu.$

\noindent Since this is true for every finite-dimensional subspace $E$ of
dimension greater than $N$ and every $\mu>\lambda_0$ the conclusion
follows.\end{proof}

\begin{Thm}\label{upper} Suppose $1\le p<\infty$.  Then there is a
Lipschitz embedding of $\ell_p$ into $c_0$ with constant $2^{1/p}$, and this
constant is best possible.
\end{Thm}

\begin{proof} The fact that $\ell_p$ does not $\lambda$-embed into $c_0$ when
$\lambda<2^{1/p}$ follows immediately from Proposition \ref{lower}. So we only
need to show that $\ell_p$ verifies condition $\Pi(2^{1/p}).$

Let $B_1$ and $B_2$ be balls with centers $a_1,a_2$ and radii $r_1,r_2.$
Suppose $\mu>2^{1/p}.$  Then $\mu<2^{1/p}(\mu^p-1)^{1/p}$. We pick $\nu$ such
that $\mu<\nu<2^{1/p}(\mu^p-1)^{1/p}$ and we fix $\epsilon>0$ so that
$$2^{1/p}\big(\mu^p(r_1+r_2)^p-(r_1+r_2+2\epsilon)^p\big)^{1/p}-
2^{1+1/p}\epsilon>\nu(r_1+r_2).$$

\noindent We first select $N\in \N$ so that
$$ \sum_{k=N+1}^{\infty}|a_1(k)|^p,\sum_{k=N+1}^{\infty}|a_2(k)|^p<\epsilon^p.$$
Let $E$ be the linear span of $\{e_1,\ldots,e_N\}$ where $(e_j)_{j=1}^\infty$
is the canonical basis of $\ell_p$. Let $P$ the canonical projection of
$\ell_p$ onto $E$, $Q=I-P$ and $R=\max(\|a_1\|+r_1,\|a_2\|+r_2).$ Then we
partition $RB_E$ into finitely many sets $A_1,\ldots,A_m$ with $\diam
A_j<\epsilon.$

\noindent Now, set $U_j=\{x\in B_1:\ Px\in A_j\}$, $V_j=\{x\in B_2:\ Px\in
A_j\}$ and
$$\mathcal S=\{(j,k)\ \exists (x,y)\in U_j\times
V_k:\ \|x-y\|>\mu(r_1+r_2)\}.$$ Thus we have
$$ \{(x,y)\in B_1\times B_2:\ \|x-y\|>\mu(r_1+r_2)\}\subset
\bigcup_{(j,k)\in\mathcal S}U_j\times V_k.$$ It remains to estimate
$\delta(U_j,V_k)$ for $(j,k)\in\mathcal S.$  Suppose $u\in U_j,v\in V_k$ and
that $x\in U_j,y\in V_k$ are such that $\|x-y\|>\mu(r_1+r_2).$   Then
$$ \|u-v\| \ge \|Pu-Pv\|\ge \|Px-Py\|-2\epsilon.$$
On the other hand
$$ r_1\ge \|x-a_1\|\ge \|Qx-Qa_1\|\ge \|Qx\|-\epsilon$$ and
$$ r_2\ge \|y-a_2\|\ge \|Qy-Qa_2\|\ge \|Qy\|-\epsilon.$$
Thus
$$ \|Qx-Qy\|\le r_1+r_2+2\epsilon.$$
Now
$$ \mu^p(r_1+r_2)^p< \|Px-Py\|^p+\|Qx-Qy\|^p\le \|Px-Py\|^p+(r_1+r_2+2\epsilon)^p.$$
Hence
$$ \|Px-Py\|^p > \mu^p(r_1+r_2)^p-(r_1+r_2+2\epsilon)^p$$ and thus
$$ 2^{1/p}\delta(U_j,V_k)\ge
2^{1/p}\big(\mu^p(r_1+r_2)^p-(r_1+r_2+2\epsilon)^p\big)^{1/p}-2^{1+1/p}\epsilon
>\nu(r_1+r_2).$$
\end{proof}

We now give a second lower bound condition in place of Proposition \ref{lower}.
We do not know whether the conclusion can be improved replacing $\lambda_0^3$
by $\lambda_0.$ If $X$ has a 1-unconditional basis, $\lambda_0^3$ can be
improved to $\lambda_0^2.$

\begin{Prop}\label{altlower}  If $X$ is a separable Banach space and $f:X\to c_0$ is a Lipschitz embedding with constant $\lambda_0$ then if $\|x\|=1$ and $(x_n)_{n=1}^{\infty}$ is a normalized weakly null sequence in $X$ we have:
\begin{equation} \label{lower100}\limsup_{n\to\infty}\|x+x_n\|\le \lambda_0^3.\end{equation}
\end{Prop}

\begin{proof} We assume that $\|x-y\|\le \|f(x)-f(y)\|\le \lambda_0\|x-y\|$ for $x,y\in X.$  Let $\mathcal U$ be a non-principal ultrafilter on the natural numbers $\mathbb N.$  We start by proving that if $x\in X$ and $(y_n)_{n=1}^{\infty},(z_n)_{n=1}^{\infty}$ are two weakly null sequences with
$\lim_{n\in\mathcal U}\|y_n\|,\lim_{n\in \mathcal U}\|z_n\|\le \|x\|$ then
\begin{equation}\label{lower103} \lambda_0^{-1}\lim_{n\in\mathcal U}\|2x+y_n+z_n\|\le \lim_{m\in\mathcal U}\lim_{n\in\mathcal U}\|2x + y_m+z_n\| \le \lambda_0\lim_{n\in\mathcal U}\|2x+y_n+z_n\|.\end{equation}
Indeed it suffices to show this under the condition $\lim_{n\in\mathcal
U}\|y_n\|=\alpha,\lim_{n\in \mathcal U}\|z_n\|=\beta$ where $\alpha,\beta\le 1$
and $\|x\|=1$.  Fix any $\epsilon>0.$ Let $f(x)=(f_j(x))_{j=1}^{\infty}.$  Then
for some $N$ we have
$$ |f_j(x)-f_j(-x)|<\epsilon \qquad j>N.$$
Thus
$$ |f_j(x+y_m)-f_j(-x-z_n)|\le \lambda_0(\|y_m\|+\|z_n\|)+\epsilon \qquad j>N.$$
Hence
$$ \lim_{m\in\mathcal U}\lim_{n\in\mathcal U}\max_{j>N}|f_j(x+y_m)-f_j(-x-z_n)|\le \lambda_0(\alpha+\beta)+\epsilon.$$
and
$$ \lim_{n\in\mathcal U}\max_{j>N}|f_j(x+y_n)-f_j(-x-z_n)|\le \lambda_0(\alpha+\beta)+\epsilon.$$

Let $\sigma_j=\lim_{n\in\mathcal U}f_j(x+y_n)$ and $\tau_j=\lim_{n\in \mathcal
U}f_j(-x-z_n).$  Then
$$ \lim_{n\in\mathcal U}|f_j(x+y_n)-f_j(-x-z_n)|=|\sigma_j-\tau_j|$$ and
$$ \lim_{m\in\mathcal U}\lim_{n\in\mathcal U}|f_j(x+y_m)-f_j(-x-z_n)|=|\sigma_j-\tau_j|.$$
Thus
\begin{align*} \lim_{n\in\mathcal U}\|2x+y_n+z_n\|&\le \lim_{n\in\mathcal U}\|f(x+y_n)-f(-x-z_n)\|\\
&\le \max(\max_{1\le j\le N}|\sigma_j-\tau_j|,\lambda_0(\alpha+\beta)+\epsilon)\\
&\le \max(\lambda_0\lim_{m\in\mathcal U}\lim_{n\in\mathcal
U}\|2x+y_m+z_n\|,\lambda_0(\alpha+\beta)+\epsilon).\end{align*} Noting that
$\epsilon>0$ is arbitrary and that
$$ \alpha+\beta\le 2 \le \lim_{m\in\mathcal U}\lim_{n\in\mathcal U}\|2x+y_m+z_n\|$$ we obtain that
$$ \lim_{n\in\mathcal U}\|2x+y_n+z_n\|\le \lambda_0\lim_{m\in\mathcal U}\lim_{n\in\mathcal U}\|2x+y_m+z_n\|.$$  The other inequality in \eqref{lower103} is similar.

Now choose $x_n=y_n=-z_n$ in \eqref{lower103}.  We obtain
$$ \lim_{m\in\mathcal U}\lim_{n\in\mathcal U}\|2x + x_m-x_n\|\le 2\lambda_0\|x\|$$ provided $(x_n)_{n=1}^{\infty}$ is weakly null and $\lim_{n\in\mathcal U}\|x_n\|\le \|x\|.$
Hence
$$ \lim_{m\in\mathcal U}\|x+\tfrac12x_m\| \le \lambda_0\|x\|.$$  This inequality can be iterated to show that
$$ \lim_{m\in\mathcal U}\lim_{n\in\mathcal U}\|x+\tfrac12x_m+\tfrac12x_n\|\le \lambda_0^2\|x\|.$$
Now assume $\|x\|=1$ and $(x_n)_{n=1}^{\infty}$ is a normalized weakly null
sequence.  Then
\begin{align*}\lim_{n\in\mathcal U}\|x+x_n\|&=\frac12\lim_{n\in\mathcal U}\|2x+x_n+x_n\|\\
&\le \frac12\lambda_0 \lim_{m\in\mathcal U}\lim_{n\in\mathcal U}\|2x+x_m+x_n\|\\
&\le \lambda_0^3.\end{align*}
\end{proof}

\begin{Thm}\label{isometric} Let $X$ be a separable Banach space.
\newline
(i) If $X$ isometrically embeds into $c_0$, then $X$ is linearly isometric to a
closed subspace of $c_0$.\newline (ii) If, for every $\epsilon>0,$ $X$
Lipschitz embeds into $c_0$ with constant at most $1+\epsilon,$ then, for every
$\epsilon>0$ there is a closed subspace $Y_{\epsilon}$ of $c_0$ with
Banach-Mazur distance $d(X,Y_{\epsilon})<1+\epsilon.$\end{Thm}

 \begin{proof} (i) is a direct consequence of the result of \cite{GodefroyKalton2003} that if a {\it separable} Banach space is isometric to a subset of a Banach space $Z$ then it is also linearly isometric to a subspace of $Z.$

 (ii)  Here we observe first that if $X$ contains a subspace isomorphic to
 $\ell_1$ then, for any $\epsilon>0$, it contains a subspace $Z_\epsilon$
 with the Banach-Mazur distance
 $d(Z_\epsilon,\ell_1)\le 1+\epsilon$ by James' distortion theorem
 \cite{James1964}. Assume
 now that $X$ can be $\lambda-$embedded into $c_0$. Thus we have that for any
 $\epsilon >0$, $\ell_1$
 can be $\lambda(1+\epsilon)-$embedded into $c_0$. Then it follows from
 Aharoni's counterexample in \cite{Aharoni1974} that $\lambda\ge 2$.

 Suppose now that $X$ does not contain any isomorphic copy of $\ell_1$.
 If $\|x\|=1$ and $(x_n)_{n=1}^{\infty}$ is any normalized weakly null sequence
 we have by Proposition \ref{altlower} that
 $$ \lim_{n\to\infty}\|x+x_n\|=1.$$
 The conclusion then follows from \cite{KaltonWerner1995} Theorem 3.5.
 \end{proof}

\section{Embeddings into $c_0^+.$}

In this section and the following we complete the already thorough study of
Lipschitz embeddings into $c_0^+$ made by Pelant in \cite{Pelant1994}.

\begin{Lem}\label{Lipfunction+}  Let $(M,d)$ be a metric space and suppose
that $A,B$ and $C$ are non-empty subsets of $M.$ Then for $\epsilon>0$, there
exists a Lipschitz function $f:M\to\mathbb R_+$ with $\Lip(f)\le 1$ such that
\newline
(i) $f(x)\le \epsilon \qquad x\in C$\newline and\newline (ii) $|f(x)-f(y)|\ge
\theta=\min\big(\delta(A,B),\max(\delta(A,C),\delta(B,C))+\epsilon\big) \qquad
x\in A,\ y\in B.$
\end{Lem}

\begin{proof} Let us suppose $\delta(A,C)\ge \delta(B,C).$   Thus
$\theta=\min(\delta(A,B),\delta(A,C)+\epsilon).$  Let us define:
$$ f(x)=\max(\theta-d(x,A),0) \qquad x\in M.$$   Then
$f(x)=\theta$ for $x\in A.$  If $x\in B$ $ f(x)\le \theta-\delta(A,B)$ so that
$f(x)=0$, while if $x\in C$ we have
$$ f(x)\le \theta-\delta(A,C)\le \epsilon.$$
\end{proof}

We may now introduce a condition analogous to $\Pi(\lambda).$  We say that
$(M,d)$ has property $\Pi_+(\lambda)$, where $\lambda>1,$ if:\newline (i)
Whenever $\mu>\lambda$ there exists $\nu>\mu$ so that if $B_1$ and $B_2$ are
two metric balls of the same radius $r$, there is a finite number of sets
$(U_j)_{j=1}^N$ and $(V_j)_{j=1}^N$ so that
$$ \lambda \delta(U_j,V_j) \ge \nu r$$ and
$$ \{(x,y)\in B_1\times B_2:\ d(x,y)>\mu r\} \subset \bigcup_{j=1}^NU_j\times V_j,$$ and
\newline
(ii) If $1<\lambda\le 2,$ there exists $1<\theta<\lambda$ and a function
$\varphi:M\to [0,\infty)$ so that
\begin{equation}\label{maxcond} |\varphi(x)-\varphi(y)|\le d(x,y)\le \theta \max(\varphi(x),\varphi(y)) \qquad x,y\in M.\end{equation}

Let us note here that condition (ii) is not required when $\lambda>2$ since
fixing any $a\in X$ the function $\varphi(x)=d(x,a)$ satisfies \eqref{maxcond}
with $\theta=2.$

We can repeat the same program for property $\Pi_+(\lambda).$

\begin{Lem}\label{everymetric+}  Every metric space has property $\Pi_+(3).$
\end{Lem}

\begin{proof}  For $\mu>3,$ let
$$ U=B_1\cap\{x:\ \exists y\in B_2,\ d(x,y)> \mu r\}$$ and
$$ V=B_2\cap\{y:\ \exists x\in B_1,\ d(x,y)>\mu r\}.$$
Then $$ \{(x,y)\in B_1\times B_2:\ d(x,y)>\mu r\} \subset U\times V.$$ Suppose
$x\in U,\ y\in V.$    Then there exists $x'\in U$ with $d(x',y)>\mu r.$  Hence
$$ d(x,y)>\mu r-d(x,x')\ge (\mu-2)r.$$  Therefore we can take $\nu=3\mu-6>\mu.$
\end{proof}

\begin{Lem}\label{everylocally+} For any $\lambda>2$, every locally compact
metric space has property $\Pi_+(\lambda)$.
\end{Lem}

The proof is immediate. Let us mention that a locally compact metric space
satisfies condition (i) for every $\lambda>1$.

\medskip We also have

\begin{Lem}\label{everycompact+} For any $\lambda>1$, every compact metric
space has property $\Pi_+(\lambda)$.
\end{Lem}

\begin{proof} Let $(K,d)$ be a compact metric space. We only have to prove condition (ii).
For $\epsilon>0$, pick a finite $\epsilon$-net $F$ of $K$ and define
$\varphi_\epsilon(x)=\max(d(x,z))_{z\in F}$. For a given $\lambda>1$,
$\varphi_\epsilon$ fulfills condition (ii) of $\Pi_+(\lambda)$ if $\epsilon$ is
small enough.
\end{proof}

\begin{Prop}\label{lambdaembedding+} Suppose $\lambda_0\ge 1$ and $M$ is a metric space which Lipschitz embeds into $c_0^+$ with constant $\lambda_0.$  Then $M$ has property $\Pi_+(\lambda)$ for all $\lambda>\lambda_0.$\end{Prop}

\begin{proof} We first consider (i) of the definition of $\Pi_+(\lambda).$ Suppose $\mu>\lambda>\lambda_0.$ Let $B_1,B_2$ be metric balls of radii $r>0$ and centers $a_1,a_2.$
Let $\Delta=\{(x,y)\in B_1\times B_2:\ d(x,y)>\mu r\}.$ Let $f:M\to c_0^+$ be
an embedding such that
$$ d(x,y)\le \|f(x)-f(y)\|\le \lambda_0 d(x,y) \qquad x,y\in M.$$
Suppose $f(x)=(f_i(x))_{i=1}^{\infty}.$  Then there exists $n$ so that
$$f_i(a_1),f_i(a_2)<(\mu-\lambda)r \qquad i\ge n+1.$$
Thus if $(x,y)\in \Delta$ we have
$$ |f_i(x)-f_i(y)|\le \max(f_i(x),f_i(y))<(\mu-\lambda)r+\lambda_0 r<d(x,y), \qquad i\ge n+1.$$  Hence
$$ d(x,y)\le \max_{1\le i\le n}|f_i(x)-f_i(y)| \qquad (x,y)\in\Delta.$$

Choose $\epsilon>0$ so that $\lambda(\mu-\epsilon)>\lambda_0\mu.$ By a
compactness argument we can find coverings $(W_k)_{k=1}^m$ of $B_1$ and
$(W'_k)_{k=1}^{m'}$ of $B_2$ such that:
$$ |f_i(x)-f_i(x')|\le \tfrac12\epsilon r \qquad x,x'\in W_k,\ 1\le i\le n,\ 1\le k\le m,$$ and
$$ |f_i(x)-f_i(x')|\le \tfrac12\epsilon r\qquad x,x'\in W'_k,\ 1\le i\le n,\ 1\le k\le m'.$$
Let
$$ \mathcal S=\{(k,k')\ 1\le k\le m,\ 1\le k'\le m':\ W_k\times W'_{k'}\cap\Delta \neq \emptyset\}$$
and define $(U_j)_{j=1}^N,\ (V_j)_{j=1}^N$ in such a way that $(U_j\times
V_j)_{j=1}^N$ is an enumeration of $(W_k\times W_{k'})_{(k,k')\in\mathcal S}.$
Then $ \Delta\subset \cup_{j=1}^NU_j\times V_j$ and the same calculations as in
the proof of Proposition \ref{lambdaembedding} give that
$$ \lambda\delta(U_j,V_j)\ge \nu r\ \ \ {\text{with}}\ \
\nu=\lambda\lambda_0^{-1}(\mu-\epsilon)>\mu.$$

If $\lambda\le 2$ we also must consider (ii).  Here we define
$\varphi(x)=\lambda_0^{-1}\|f(x)\|$ where $f:M\to c_0^+$ is as above. Then
$\varphi$ satisfies \eqref{maxcond} with $\theta=\lambda_0$. Indeed,
$$|\varphi(x)-\varphi(y)|\le \lambda_0^{-1}\|f(x)-f(y)\|\le d(x,y)$$
and
$$d(x,y)\le \|f(x)-f(y)\| \le \max(\|f(x)\|,\|f(y)\|) \le
\lambda_0\max(\varphi(x),\varphi(y)).$$

\end{proof}

\medskip Next, in place of Lemma \ref{soupup} we have

\begin{Lem}\label{soupup+}  Let $\lambda>1$ and $(M,d)$ be a metric space with
property $\Pi_+(\lambda).$ Then for every $\mu>\lambda$ there is a constant
$\nu>\mu$ so that if  $B_1$ and $B_2$ are two metric balls of radius $r$ then
there are finitely many sets $(U_j)_{j=1}^N,(V_j)_{j=1}^N$ such that if
$(x,y)\in B_1\times B_2$ and $d(x,y)>\mu r$ then there exists $1\le j\le N$ so
that $x\in U_j,\ y\in V_j$ and:
$$  \lambda \mu\delta(U_j,V_j) \ge \nu d(x,y).$$ \end{Lem}

We omit the proof of this which is very similar to that of Lemma \ref{soupup}
and only uses part (i) of the definition of $\Pi_+(\lambda)$.

\medskip Then we have the following analogue of Lemma \ref{construction}.

\begin{Lem}\label{construction+} Let $\lambda >1$. Suppose $(M,d)$ has property
$\Pi_+(\lambda)$. Suppose $0<\alpha<\beta$.  Let $F,G$ be finite subsets of $M$
and let $\Delta_+(F,G,\alpha,\beta)$ be the set of $(x,y)\in M\times M$ such
that
$$ \lambda \max(d(x,G),d(y,G))+\alpha\le d(x,y)<\lambda\max(d(x,F),d(y,F))+\beta.$$
Then there is a finite set $\mathcal F=\mathcal F(F,G,\alpha,\beta)$ of
functions $f:M\to\mathbb R_+$ with $\Lip(f)\le \lambda$ and such that
$$ f(x)\le \lambda\beta \qquad x\in F$$ and
$$ d(x,y)< \max_{f\in \mathcal F}|f(x)-f(y)| \qquad (x,y)\in \Delta_+(F,G,\alpha,\beta).$$
\end{Lem}

\begin{proof} We first argue that for some constant $K$ we have
$$ d(x,y)\le K,\qquad x,y\in \Delta_+(F,G,\alpha,\beta).$$
If $\lambda>2$ this is follows from the fact that
$$d(x,G)+d(y,G)\ge d(x,y)-R$$ where $R$ is the diameter of $G$. Hence
$$d(x,y)\leq K=\lambda(\lambda-2)^{-1}R,\qquad x,y\in
\Delta_+(F,G,\alpha,\beta).$$ In the case $1< \lambda\le 2$ let
$\varphi,\theta$ be as in the definition of $\Pi_+(\lambda)$ and satisfy
\eqref{maxcond}.  Let $K_0=\max\{\varphi(z):\ z\in G\}.$ Thus
\begin{align*} \lambda d(x,y)&\le \lambda\theta \max(\varphi(x),\varphi(y))\\
&\le \lambda \theta K_0+\lambda\theta\max(d(x,G),d(y,G))\\
&\le \lambda \theta K_0 +\theta d(x,y) \qquad x,y\in
\Delta_+(F,G,\alpha,\beta),\end{align*} so that
$$ d(x,y) \le K=\frac{\lambda \theta K_0}{\lambda-\theta},
\qquad x,y\in \Delta_+(F,G,\alpha,\beta).$$

We next let
$$ \mu=\lambda +\frac{\alpha\lambda}{2K}$$ and choose
$\nu=\nu(\mu)$ according to the conclusion of Lemma \ref{soupup+}. We fix
$\epsilon>0$ so that
$\epsilon<\min(\frac{\alpha}{2},\lambda^{-1}(\lambda-1)\beta).$

Let $E=\{x:\ d(x,G)\le \lambda^{-1}K\}.$  Since $E$ is metrically bounded and
$F\cup G$ is finite we can partition $E$ into finitely many subsets
$(E_1,\ldots,E_m)$ so that for each $z\in F\cup G$ we have:
$$ |d(x,z)-d(x',z)|\le \epsilon \qquad x,x'\in E_j,\ 1\le j\le m.$$

\noindent For each $j$, we define $z_j\in G$ and $r_j$, as in the proof of
Lemma \ref{construction}, so that

$$ \inf_{x\in E_j}d(x,z_j)=\inf_{x\in E_j}d(x,G)=r_j.$$
Note that $r_j\le \lambda^{-1}K$ and $E_j$ is contained in a ball $B_j$
centered at $z_j$ with radius $r_j+\epsilon.$

\noindent Now for each pair $(j,k)$ we denote $B_{j,k}$ the ball with center
$z_j$ and radius $\max(r_j+\epsilon,r_k+\epsilon)$. By Lemma \ref{soupup+}, we
can find finitely many pairs of sets $(\tilde U_{jkl},\tilde
V_{jkl})_{l=1}^{N_{jk}}$ such that for every $(x,y)\in B_{j,k}\times B_{k,j}$
with $ d(x,y)>\mu(\max(r_j,r_k)+\epsilon)$ there exists $1\le l\le N_{jk}$ with
$x\in \tilde U_{jkl},y\in \tilde V_{jkl}$ and
$$ \lambda\mu\delta(\tilde U_{jkl},\tilde V_{jkl})\ge \nu d(x,y).$$
Then we set $U_{jkl}=\tilde U_{jkl} \cap E_j$ and $V_{jkl}=\tilde V_{jkl}\cap
E_k$.

We now apply Lemma \ref{Lipfunction+} to construct Lipschitz functions
$f_{jkl}:M\to\mathbb R_+$ where $1\le j,k\le m,\ 1\le l\le N_{jk}$ such that
$\Lip(f_{jkl})\le \lambda,$
$$ f_{jkl}(x)\le \lambda \beta \qquad x\in F$$ and
$$ |f_{jkl}(x)-f_{jkl}(y)|\ge \lambda\theta_{jkl} \qquad x\in U_{jkl},\ y\in V_{jkl}$$ where
$$ \theta_{jkl}=\min\big(\delta(U_{jkl},V_{jkl}),\max(\delta(U_{jkl},F),\delta(V_{jkl},F))+\beta\big).$$

Now let us suppose $(x,y)\in \Delta_+(F,G,\alpha,\beta).$  Then there exists
$(j,k)$ so that $x\in E_j,\ y\in E_k.$  It follows from our choice of $\mu$ and
$\epsilon$ that
\begin{align*} d(x,y)&\ge \lambda\max(d(x,G),d(y,G))+\alpha\\
&\ge \lambda \max(r_j,r_k)+\alpha>\mu\max(r_j,r_k)+\epsilon.\end{align*}

\noindent Thus  there exists $1\le l\le N_{jk}$ so that $x\in U_{jkl},y\in
V_{jkl}$ and
$$ \lambda\delta(U_{jkl},V_{jkl})\ge \frac{\nu}{\mu} d(x,y)>d(x,y).$$
On the other hand, $\epsilon<\lambda^{-1}(\lambda-1)\beta$, so
\begin{align*} \lambda\max(\delta(U_{jkl},F),\delta(V_{jkl},F))+\beta)
&\ge \lambda\max(d(x,F),d(y,F))+\lambda(\beta-\epsilon)\\
&>\lambda \max(d(x,F),d(y,F))+\beta\\ &> d(x,y).\end{align*} Hence
$$ |f_{jkl}(x)-f_{jkl}(y)|\ge \lambda \theta_{jkl}>d(x,y).$$  Thus we can take for $\mathcal F$
the collection of all functions $f_{jkl}$ for $1\le j,k\le m,\ 1\le l\le
N_{jk}.$\end{proof}

Finally our theorem is

\begin{Thm}\label{Pi+} Suppose a separable metric space $(M,d)$ has property $\Pi_+(\lambda).$
Then there is a Lipschitz embedding $f:M\to c_0^+$ with
$$ d(x,y)<\|f(x)-f(y)\|\le \lambda d(x,y) \qquad x,y\in M,\ x\neq y.$$
\end{Thm}

\begin{proof} We use the notation of the proof of Theorem \ref{Pi}.
Then we build an increasing sequence of integers $(n_k)_{k=0}^{\infty}$ (with
$n_0=0$) and a sequence $(f_j)_{j=1}^{\infty}$ of Lipschitz functions
$f_j:M\to\mathbb R_+$ with $\Lip(f_j)\le \lambda$ so that
$$ f_j(x)\le \lambda\epsilon_k \qquad x\in F_k, \ n_{k-1}<j\le n_k$$ and
if
\begin{equation}\label{selectk+} \lambda \max(d(x,F_{k+1}),d(y,F_{k+1}))+
\epsilon_{k+1}\le d(x,y)<\lambda \max(d(x,F_k),d(y,F_k))+\epsilon_k
\end{equation}
then
$$ \max_{n_{k-1}<j\le n_k}|f_j(x)-f_j(y)|>d(x,y).$$
If $x\neq y$ the sequence
$$ \tau_k=\lambda\max(d(x,F_k),d(y,F_k))+\epsilon_k$$ is decreasing
and tends to zero.

\noindent If $\lambda>2$, we clearly have $\tau_1>d(x,y)$.

\noindent Assume $1<\lambda\le 2$. Let $\varphi$ be given by the part (ii) of
property $\Pi_+(\lambda)$. We choose $\epsilon_1>\lambda\varphi(u_1)$. Then we
have
$$d(x,y)\le \lambda\max(\varphi(x),\varphi(y))<
\epsilon_1+\lambda\max(d(x,u_1),d(y,u_1))=\tau_1.$$ Hence, in both cases the
desired embedding can be defined again by $f(x)=(f_j(x))_{j=1}^{\infty}.$
\end{proof}

As a first corollary, we obtain the two following results due to Pelant
(\cite{Pelant1994}).

\begin{Cor}\label{Pelant}

(a) For every separable metric space $(M,d)$ there is a Lipschitz embedding
$f:M\to c_0^+$ so that
$$ d(x,y)<\|f(x)-f(y)\|\le 3 d(x,y) \qquad x,y\in M,\ x\neq y.$$

\noindent (b) For any compact metric space $(K,d)$ and any $\lambda>1$, $(K,d)$
$\lambda$-embeds into $c_0^+$.
\end{Cor}

It is proved in \cite{Pelant1994} that both of the above statements are
optimal. This was also known to Aharoni \cite{Aharoni1978} for part (a).

\medskip We also have.

\begin{Thm}\label{embedlocal+} For every locally compact metric space $(M,d)$
and every $\lambda>2$, $(M,d)$ $\lambda$-embeds into $c_0^+$. This result is
optimal.
\end{Thm}
\begin{proof} The result is obtained by combining Theorem \ref{Pi+} and Lemma
\ref{everylocally+}. We only have to show its optimality.

Let $\cal D$ be the set of all finite sequences with values in $\{0,1\}$
including the empty sequence denoted $\emptyset$ and let $\cal D^*=\cal
D\setminus \{\emptyset\}$. For $s\in \cal D$, we denote $|s|$ its length. Then
$(e_s)_{s\in \cal D}$ is the canonical basis of $\ell_1(\cal D)$. We consider
the following metric subspace of $\ell_1(\cal D)$:
$$M=\{0,e_\emptyset\}\cup\{|s|e_s,e_\emptyset+|s|e_s,\ s\in \cal D^*\}.$$
This is clearly a locally finite metric space. Assume now that there exists
$f=(f_k)_{k=1}^\infty:M \to c_0^+$ such that
$$\|x-y\|_1\le \|f(x)-f(y)\|_\infty \le 2\|x-y\|_1 \qquad  x,y \in M.$$
There exits $K\ge 1$ such that $f_k(e_\emptyset)<1$ and $f_k(0)<1$ for all
$k>K$. Then, using the positivity of $f$, we obtain
\begin{align*} |f_k(e_\emptyset+ne_s)-f_k(ne_t)|&\le
\max\big(f_k(e_\emptyset+ne_s),f_k(ne_t)\big)\\& < 1+2n  \qquad k>K,\ s\neq t,\
|s|=|t|=n.\end{align*} On the other hand,
$$\|f(e_\emptyset+ne_s)-f(ne_t)\|_\infty \ge 1+2n \qquad s\neq t,\
|s|=|t|=n.$$ Thus, for all $s\neq t,\ |s|=|t|=n$, there exists $k\le K$ so that
$$|f_k(e_\emptyset+ne_s)-f_k(ne_t)|\ge 1+2n.$$
Let now $C=\max(\|f(e_\emptyset)\|_\infty,\|f(0)\|_\infty)$. Then
$$|f_k(e_\emptyset+ne_s)-f_k(ne_t)|\le C+2n  \qquad k\le K,\ s\neq t,\ |s|=|t|=n.$$
Thus, for $n$ large enough and all $s\neq t,\ |s|=|t|=n$, there exists $k\leq
K$ such that either
$$f_k(ne_s)\le C-1 \qquad {\rm and} \qquad f_k(e_\emptyset+ne_t)\ge 1+2n$$
or
$$f_k(ne_s)\ge 1+2n \qquad {\rm and} \qquad f_k(e_\emptyset+ne_t)\le C-1.$$
Therefore: either
$$f_k(ne_s)\le C-1 \qquad {\rm and} \qquad f_k(ne_t)\ge 2n-1$$
or
$$f_k(ne_s)\ge 1+2n \qquad {\rm and} \qquad f_k(ne_t)\le C+1.$$
Let us now denote $\alpha(k,s)=\car_{[0,C+1]}(f_k(|s|e_s))$. Then, for $n$ big
enough, we have that for all $s\neq t,\ |s|=|t|=n$, there exists $k\le K$ so
that $\alpha(k,s)\neq\alpha(k,t)$. This is clearly impossible if $n>K$. This
finishes our proof.

\end{proof}

\section{Embeddings of subsets of classical Banach spaces into $c_0^+$.}

\begin{Prop}\label{lower+} Suppose $X$ is a separable Banach space and that $f:X\to c_0^+$ is a
Lipschitz embedding with constant $\lambda_0.$  Then for any $u\in X$ with
$\|u\|=1$ and any infinite-dimensional subspace $Y$ of $X$ we have
$$ \inf_{\substack{y\in Y\\ \|y\|=1}}\|u+ 2y\|\le \lambda_0.$$
\end{Prop}

\begin{proof} The proof is almost identical to that of Proposition \ref{lower}. It follows from Proposition \ref{lambdaembedding+} that $X$ has
property $\Pi_+(\lambda)$ for any $\lambda>\lambda_0$. We consider
$B_1=-u+2B_X$ and $B_2=u+2B_X$, where $B_X$ denotes the closed unit ball of
$X$. Suppose $\mu>\lambda_0$ and select $\mu>\lambda>\lambda_0.$ Then, for some
$\nu>\mu$, we can find finitely many closed sets $(U_j,V_j)_{j=1}^N$ verifying:
$$\lambda \delta(U_j,V_j)\ge 2\nu$$ and
$$ \{(x,y)\in B_1\times B_2:\ \|x-y\|>2\mu\}\subset \bigcup_{j=1}^NU_j\times V_j.$$
Now let $E$ be any subspace of $X$ of dimension greater than $N$ and let
$$ A_j=\{e\in E:\ \|e\|=1,\ (-u+2e,u-2e)\in {U}_j\times {V}_j\}.$$  We then conclude the proof as in
Proposition \ref{lower}. Assume that for any $e\in S_E$, $\|u+2e\|>\mu$. Then
$A_1\cup\cdots\cup A_N=S_E$ and so there exists $e$ in $S_E$ and $k\leq N$ such
that $e$ and $-e$ belong to $A_k$, i.e. $-u\pm 2e\in U_k$ and $u\pm 2e\in V_k.$
This implies that $\delta (U_k,V_k)\le 2$, which is a contradiction. So, there
exists $e\in S_E$ with $\|u+ 2e\|\le \mu$ and we conclude as in the proof of
Proposition \ref{lower}.
\end{proof}

\begin{Thm}\label{ellp+}  Suppose $1\le p<\infty$.\newline
(i) There is a Lipschitz embedding of $\ell_p$ into $c_0^+$ with constant
$(2^p+1)^{1/p}$ and this is best possible.\newline (ii) There is a Lipschitz
embedding of $\ell_p^+$ into $c_0^+$ with constant $3^{1/p}$ and this is best
possible.
\end{Thm}

\begin{proof}  Let us prove first that $\ell_p$ has $\Pi_+(c_p)$ where $c_p=(1+2^p)^{1/p}.$  The proof is very similar to that of Theorem \ref{upper}.
Let $B_1$ and $B_2$ be balls with centers $a_1,a_2$ and radius $r>0.$ Suppose
$\mu>c_p$ and that suppose $\mu<\nu<c_p(\mu^p-2^p)^{1/p}.$ Fix $\epsilon>0$
such that
$$c_p\big(\mu^p r^p-2^p(r+\epsilon)^p\big)^{1/p}-2\epsilon c_p>\nu r.$$

\noindent We select $N\in \N$ so that
$$ \sum_{k=N+1}^{\infty}|a_1(k)|^p,\sum_{k=N+1}^{\infty}|a_2(k)|^p<\epsilon^p.$$
Let $E$ be the linear span of $\{e_1,\ldots,e_N\}$ where $(e_j)$ is the
canonical basis of $\ell_p$. Let $P$ the canonical projection of $\ell_p$ onto
$E$, $Q=I-P$ and $R=\max(\|a_1\|,\|a_2\|)+r.$ Then we partition $RB_E$ into
finitely many sets $A_1,\ldots,A_m$ with $\diam A_j<\epsilon.$

\noindent Now, set $U_j=\{x\in B_1:\ Px\in A_j\}$, $V_j=\{x\in B_2:\ Px\in
A_j\}$ and
$$\mathcal S=\{(j,k):\ \exists (x,y)\in U_j\times
V_k:\ \|x-y\|>\mu r\}.$$ Thus we have
$$ \{(x,y)\in B_1\times B_2:\ \|x-y\|>\mu r\}\subset
\bigcup_{(j,k)\in\mathcal S}U_j\times V_k.$$ It remains to estimate
$\delta(U_j,V_k)$ for $(j,k)\in\mathcal S.$  Suppose $u\in U_j,v\in V_k$ and
that $x\in U_j,y\in V_k$ are such that $\|x-y\|>\mu r.$ Then
$$ \|u-v\| \ge \|Pu-Pv\|\ge \|Px-Py\|-2\epsilon.$$
On the other hand
$$ r\ge \|x-a_1\|\ge  \|Qx\|-\epsilon$$ and
$$ r\ge \|y-a_2\|\ge  \|Qy\|-\epsilon.$$
Thus
\begin{equation}\label{equn10} \|Qx-Qy\|\le 2r+2\epsilon.\end{equation} Now
$$ \mu^p r^p< \|Px-Py\|^p+\|Qx-Qy\|^p\le \|Px-Py\|^p+ 2^p(r+\epsilon)^p.$$
Hence
$$ \|Px-Py\|^p > \mu^p r^p-2^p(r+\epsilon)^p,$$ and so
$$ c_p\delta(U_j,V_k)\ge
c_p\big(\mu^p r^p-2^p(r+\epsilon)^p\big)^{1/p}-2\epsilon c_p>\nu r.$$ Hence
$\ell_p$ has $\Pi_+(c_p).$

Next we show that $\ell_p^+$ has property $\Pi_+(3^{1/p})$.  To do this we
repeat the argument above. We take $\mu>3^{1/p}$ and suppose that
$\mu<\nu<3^{1/p}(\mu^p-2)^{1/p}.$  Choose $\epsilon>0$ so that:
$$3^{1/p}\big(\mu^pr^p-2(r+\epsilon)^p\big)-2\epsilon3^{1/p}>\nu r.$$
Next repeat the construction, but working inside the positive cone $\ell_p^+$.
The only difference is that \eqref{equn10} is replaced by
\begin{equation}\label{equn11} \|Qx-Qy\|\le
2^{1/p}\max(\|Qx\|,\|Qy\|)\le 2^{1/p}(r+\epsilon).\end{equation} Hence
$$ \|Px-Py\|^p>\mu^p r^p-2(r+\epsilon)^p,$$ and so this time
$$ 3^{1/p}\delta(U_j,V_k) \ge 3^{1/p}\big(\mu^pr^p-2(r+\epsilon)^p\big)-
2\epsilon3^{1/p}>\nu r.$$

For the second half of the condition when $3^{1/p}\le 2$ we note that
$\varphi(x)=\|x\|$ satisfies \eqref{maxcond} with $\theta=2^{1/p}<3^{1/p}.$

These calculations combined with Theorem \ref{Pi+} show the existence of the
Lipschitz embeddings in parts (i) and (ii). Proposition \ref{lower+} shows the
constant is best possible when in (i).  For (ii) let us suppose $f:\ell_p^+\to
c_0^+$ is an embedding such that
$$ \|x-y\|\le \|f(x)-f(y)\|\le \lambda\|x-y\| \qquad x,y\in \ell_p^+$$ where $\lambda<3^{1/p}.$
Let $f(x)=(f_j(x))_{j=1}^{\infty}.$  Let $\epsilon=(3^{1/p}-\lambda)/2.$  Then
there exists $N$ such that
$$\max(f_j(e_1),f_j(0))<\epsilon \qquad j\ge N+1.$$
Hence if $m,n>1$
$$ |f_j(e_1+e_m)-f_j(e_n)|\le \max (f_j(e_1+e_m),f_j(e_n))\le \lambda +\epsilon<3^{1/p}, \quad j\ge N+1.$$
Now we may pass to a subsequence so that the following limits exist:
$$ \lim_{k\to\infty}f_j(e_1+e_{n_k})=\sigma_j,\quad \lim_{k\to\infty}f_j(e_{n_k})=\tau_j,\qquad 1\le j\le N.$$
Clearly
$$ |\sigma_j-\tau_j|\le \lambda, \qquad 1\le j\le N.$$
Now
$$ \lim_{k\to\infty}|f_j(e_1+e_{n_k})-f_j(e_{n_{k+1}})| \le \lambda \qquad 1\le j\le N$$ and we have a contradiction since $\|e_1+e_{n_k}-e_{n_{k+1}}\|=3^{1/p}>\lambda.$
\end{proof}

\section{Spaces embedding isometrically into $c_0$ and $c_0^+$.}

In this final section we study isometric embeddings into $c_0$ and $c_0^+.$
Note that a separable Banach space isometrically embeds into $c_0$ if and only
if it embeds linearly and isometrically \cite{GodefroyKalton2003}.

We recall that a metric space $(M,d)$ is an ultrametric space if
\begin{equation} \label{ultrametric} d(x,y)\le \max(d(x,z),d(z,y))
\qquad x,y,z\in M.\end{equation} Note that this implies
\begin{equation} \label{ultrametric1} d(x,y)=\max(d(x,z),d(z,y)) \qquad d(x,z)\neq d(z,y).\end{equation}

\begin{Lem}\label{Gamma} Let $(M,d)$ be a separable ultrametric space.  Then there is a countable subset $\Gamma$ of $[0,\infty)$ such that $d(x,y)\in\Gamma$ for all $x,y\in M.$\end{Lem}

\begin{proof} For each fixed $x\in M$ let
$\Gamma_x=\{d(x,y):y\in M\}.$  Suppose $\Gamma_x$ is uncountable; then for some
$\delta>0$ the set $\Gamma_x\cap (\delta,\infty)$ is uncountable. Pick an
uncountable set $(y_i)_{i\in I}$ in $M$ so that $d(x,y_i)>\delta$ and the
values of $d(x,y_i)$ are distinct for $i\in I.$  Then $i\neq j$ we have
$d(y_i,y_j)>\delta$ by \eqref{ultrametric1}.  This contradicts separability of
$M.$

Thus each $\Gamma_x$ is countable.  Let $D$ be a countable dense subset of $M$
and let $\Gamma=\cup_{x\in D}\Gamma_x.$  If $y,z\in M$ with $y\neq z,$ pick
$x\in D$ with $d(x,y)<d(y,z).$  Then $d(y,z)=d(x,z)\in\Gamma$ by
\eqref{ultrametric1}.\end{proof}

\begin{Thm}\label{ultrametric2}
Every separable ultrametric space embeds isometrically into $c_0^+$\end{Thm}

\begin{proof} Pick $\Gamma$ as in Lemma \ref{Gamma}.  Let
$(a_j)_{j=1}^{\infty}$ be a countable dense subset of an ultrametric space $M$.
Let $\mathcal D$ be the collection of finite sequences $(r_1,\ldots,r_n)$ with
$r_j\in\Gamma$ for $1\le j\le n.$ For each $(r_1,\ldots,r_n)\in\mathcal D$ we
define a function $f_{r_1,\ldots,r_n}$ by
$$ f_{r_1,\ldots,r_n}(x)=\begin{cases}\min(r_1,\ldots,r_n) \qquad d(x,a_j)=r_j,\ 1\le j\le n\\
0 \qquad \qquad \qquad \qquad \text{otherwise}.\end{cases}$$ If $x\in M$ let
$d(x,a_j)=s_j.$  Then $\lim_{n\to\infty}\min(s_1,\ldots,s_n)=0$ and it follows
that $f(x)=(f_{r_1,\ldots,r_n}(x))_{(r_1,\ldots,r_n)\in\mathcal D}$ is a map
from $M$ into $c_0^+(\mathcal D).$

If $x,y\in M$ and $f_{r_1,\ldots,r_n}(x)\neq f_{r_1,\ldots,r_n}(y)$ we can
assume without loss of generality that $d(x,a_j)=r_j$ for $1\le j\le n$ but
that for some $1\le k\le n$ we have $d(y,a_k)\neq r_k.$ Then
$$ |f_{r_1,\ldots,r_n}(x)-f_{r_1,\ldots,r_n}(y)|=\min(r_1,\ldots,r_n)
\le r_k\le \max(d(x,a_k),d(y,a_k))= d(x,y)$$ by \eqref{ultrametric1}. Thus
$\|f(x)-f(y)\|\le d(x,y)$ for $x,y\in M.$

On the other hand if $x\neq y$ there is a least $k$ so that $d(x,a_k)\neq
d(y,a_k).$  Assume $d(x,a_k)>d(y,a_k)$ and $r_j=d(x,a_j)$ for $1\le j\le k.$
Then $ d(x,y)=r_k.$  On the other hand $d(x,y)\le r_j$ for $1\le j\le k.$ Hence
$$d(x,y)=r_k=|f_{r_1,\ldots,r_k}(x)-f_{r_1,\ldots,r_k}(y)|.$$  Thus
$f$ is an isometry.
\end{proof}

As a final example we consider an infinite branching tree $\mathcal T$ defined
as the set of all ordered subsets (nodes) $a=(m_1,\ldots,m_k)$ (where
$m_1<m_2<\ldots<m_k$) of $\mathbb N$ (including the empty set).  Let $|a|=k$ be
the length of $a$ so that $|\emptyset|=0.$  If
$a=(m_1,\ldots,m_k),b=(n_1,\ldots,n_l)$ are two nodes we define $a\wedge b$ to
be the node $(m_1,\ldots,m_r)$ where $r\le \min(k,l)$ is the greatest integer
such that $m_j=n_j$ for $1\le j\le r.$ We write $a\prec b$ if $b\wedge a=a.$
$\mathcal T$ is a graph if we define two nodes $a,b$ to be adjacent if
$||a|-|b||=1$ and $a\prec b$ or $b\prec a$.  The natural graph metric $d$ is
thus given by
$$ d(a,b)=|a|+|b|-2|a\wedge b|.$$

\begin{Thm}  The infinite branching tree embeds isometrically into $c_0.$
\end{Thm}

\begin{proof} For each $(a,n)\in\mathcal T\times \mathbb N$ we define
$$ f_{a,n}(b)=\begin{cases} |b|-|a| \qquad a\prec b,\ b\neq a,\ b_{|a|+1}=n\\
|a|-|b| \qquad a\prec b,\ b\neq a,\ b_{|a|+1}>n\\
0 \qquad\qquad \text{otherwise}.\end{cases}$$ For fixed $b$ we have
$f_{a,n}(b)\neq 0$ only when $a\prec b$ and $n\le b_{|a|+1}$ and this is a
finite set.  Hence $f(b)=(f_{a,n}(b))_{(a,n)\in\mathcal T\times\mathbb N}$
defines a map of $\mathcal T$ into $c_0(\mathcal T\times \mathbb N).$

Suppose $d(b,b')=1$ and that $|b'|=|b|+1.$  Then by examining cases it is clear
that $|f_{a,n}(b)-f_{a,n}(b')|\le 1$ so that $\|f(b)-f(b')\|\le 1.$  It follows
that $\|f(b)-f(b')\|\le d(b,b')$ for arbitrary $b,b'\in\mathcal T.$

\noindent If $b\neq b'$ pick $a=b\wedge b'$ and assume as we may that either
that $b'=a\wedge b=a$ or $b_{|a|+1}<b'_{|a|+1}.$  Put $n=b_{|a|+1}$. Then
$$ f_{a,n}(b)=|b|-|a|,\quad f_{a,n}(b')=|a|-|b'|$$ so that
$$ |f_{a,n}(b)-f_{a,n}(b')|=d(b,b').$$  Hence $f$ is an isometry. \end{proof}

\noindent {\it Remark}. Since $c_0$ 2-embeds into $c_0^+$, so does $\mathcal
T$. It follows from the fact that $\mathcal T$ contains a copy of $\mathbb Z$,
that it is again optimal.


\begin{bibsection}
\begin{biblist}
\bib{Aharoni1974}{article}{
  author={Aharoni, I.},
  title={Every separable metric space is Lipschitz equivalent to a subset of $c\sp {+}\sb {0}$},
  journal={Israel J. Math.},
  volume={19},
  date={1974},
  pages={284--291},
}

\bib{Aharoni1978}{book}{
  author={Aharoni, I.},
  title={Lipschitz maps and uniformly continuous functions between Banach spaces},
  series={Ph.D. thesis},
  place={Hebrew University, Jerusalem},
  date={1978},
}

\bib{Assouad1978}{article}{
  author={Assouad, P.},
  title={Remarques sur un article de Israel Aharoni sur les prolongements lipschitziens dans $c\sb {0}$ (Israel J. Math. 19 (1974), 284--291)},
  journal={Israel J. Math.},
  volume={31},
  date={1978},
  pages={97--100},
}

\bib{BenyaminiLindenstrauss2000}{book}{
  author={Benyamini, Y.},
  author={Lindenstrauss, J.},
  title={Geometric nonlinear functional analysis. Vol. 1},
  series={American Mathematical Society Colloquium Publications},
  volume={48},
  publisher={American Mathematical Society},
  place={Providence, RI},
  date={2000},
}

\bib{GodefroyKalton2003}{article}{
  author={Godefroy, G.},
  author={Kalton, N. J.},
  title={Lipschitz-free Banach spaces},
  journal={Studia Math.},
  volume={159},
  date={2003},
  pages={121\ndash 141},
}

\bib{James1964}{article}{
  author={James, R. C.},
  title={Uniformly non-square Banach spaces},
  journal={Ann. of Math. (2)},
  volume={80},
  date={1964},
  pages={542\ndash 550},
}

\bib{KaltonWerner1995}{article}{
  author={Kalton, N. J.},
  author={Werner, D.},
  title={Property $(M)$, $M$-ideals, and almost isometric structure of Banach spaces},
  journal={J. Reine Angew. Math.},
  volume={461},
  date={1995},
  pages={137\ndash 178},
}

\bib{LyusternikShnirelman1930}{book}{
  author={Lyusternik, L.},
  author={Shnirelman, S.},
  title={Toplogical methods in variational problems (Russian)},
  publisher={Iss. Institut Matem. i Mech. pri O.M.G.U.},
  place={Moscow},
  date={1930},
}

\bib{Matousek2003}{book}{
  author={Matou{\v {s}}ek, J.},
  title={Using the Borsuk-Ulam theorem},
  series={Universitext},
  publisher={Springer-Verlag},
  place={Berlin},
  date={2003},
  pages={xii+196},
}

\bib{Pelant1994}{article}{
  author={Pelant, J.},
  title={Embeddings into $c\sb 0$},
  journal={Topology Appl.},
  volume={57},
  date={1994},
  pages={259--269},
}


%\bibselect{references2}
\end{biblist}
\end{bibsection}
\end{document}